\documentclass[a4paper]{amsart}

 \title[Analytic version of Hilbert's $17$-th problem]{A non-commutative,\\
   analytic version of Hilbert's  $17$-th problem\\ in type II$_1$ 
von Neumann algebras\footnote{Research partially supported by NSF GRANT no. DMS 0200741}}
 \author{Florin R\u{a}dulescu}

\address{Department of Mathematics \\
University of Roma ``Tor Vergata''\\ Via della Ricerca Scientifica, 00133
Roma, Italy\\
on leave from University of Iowa \\ Iowa City, IA 52242, USA}

 \usepackage{amsmath,amsthm}
 \usepackage{amsmath}
 \usepackage{amsfonts}
 \usepackage{amssymb}
 \usepackage{amsthm}
 \usepackage{calc}
 \usepackage{times}

 \def\prep{\perp}
 \def\trr{\mathop{\rm Tr}}
\def\tr{\mathop{\rm tr}}

 \def\leq{\leqslant}
 \def\geq{\geqslant}

 \def\epsilon{\varepsilon}

\def\cI{{\cal I}}
\def\cH{{\cal H}}
 \def\C{{\mathbb C}}
 \def\N{{\mathbb N}}
 \def\R{{\mathbb R}}
 \def\diez{\sharp}

 \theoremstyle{plain} 
 \newtheorem{thm}{Theorem}[section]
 \newtheorem{cor}[thm]{Corollary}
 \newtheorem{lem}[thm]{Lemma}
 \newtheorem{prop}[thm]{Proposition}

 \theoremstyle{definition}

 \newtheorem{obs}[thm]{Observation}

 \def\can{{\mathbb C}_{\rm an}}
 \def\cal{\mathcal}

 \subjclass[2000]{Primary: 46L05; Secondary: 46L10}


\begin{document}

 \begin{abstract}
 Let $Y_1, \ldots, Y_n$ be  $n$ indeterminates. For $I=(i_1, \ldots,i_p)$,
 $i_s\in \{1, 2$, $\ldots, n\}$, $s= 1, 2, \ldots, p$, let $Y_I$ be the
 monomial $Y_{i_1} Y_{i_2} \cdots Y_{i_p}$. Denote by $|I|=p$. Let $\can[Y_1, Y_2, \ldots, Y_p]$ be the ring of non-commutative  series $\sum a_I Y_I$, $a_I\in \C$, such that 
 $\sum |a_I| R^{|I|}<\infty$ for all $R>0$. On $\can[Y_1, Y_2, \ldots, Y_n]$ we have a canonical involution extending by linearity $(a_I Y_I)^* = \overline a Y_{I^{\rm op}}$, $a_I\in \C$, $I\in \cI_n$, $I=\{i_1, i_2,  \ldots,i_p\}$, $I^{\rm op}= \{i_p, i_{p-1}. \ldots, i_1\}$. By $\can^{\rm sym}[Y_1, Y_2, \ldots,Y_n]$ we denote the real subspace of $\can [Y_1, Y_2, \ldots, Y_n]$ of series that are auto-adjoint.
 We say that two series $p, q$ are cyclic equivalent if $p-q$ is a sum
 (possible infinite) of scalar multiples of monomials of the type
 $Y_I-Y_{\tilde I}$, where $\widetilde I$ is a cyclic permutation of $I$. We
 call a series $q$ in $\can [Y_1, \ldots, Y_p]$ a sum of squares if $q$
 is a weak limit of sums $\sum_s b^*_s b_s$, where $b_s\in \can [Y_1, \ldots, Y_p]$.

 We prove that if a series $p(Y_1, \ldots,Y_n)$ in $\can^{\rm sym}[Y_1, \ldots,
 Y_n]$ has the property that $\tau (p(X_1, \ldots, X_n))\geq 0$ for every $M$
 type II$_1$ von Neumann algebra with faithful trace $\tau$ and for all selfadjoint
 $X_1, X_2, \ldots, X_n$ in $M$, then $p$ is equivalent to a sum of squares in
 $\can[Y_1, \ldots, Y_n]$. As a corollary, it follows that the Connes embedding
 conjecture is equivalent to a statement on the structure of matrix trace
 inequalities: if $p(Y_1, \ldots, Y_n)$ in $\can^{\rm sym} [Y_1, \ldots, Y_n]$
 is such that $\tr p(X_1, \ldots,X_n)\geq 0$, for all selfadjont matrices $X_1,
 \ldots, X_n$, of any size, then $p$ should be equivalent to a sum of squares
  in\linebreak $\can [Y_1, \ldots,Y_n]$.  
 \end{abstract}

 \maketitle

 \thispagestyle{empty}

 \section{Introduction}
 The Connes embedding conjecture (\cite{Co}) states that:
 \vskip4pt

  {\it Every type {\rm II}$_1$
   factor (equiva\-lently any type {\rm II}$_1$ von Neumann algebra) can be
   embedded into the factor $R^\omega$} (see \cite{Co}, p.~105).
 \vskip4pt

 Equivalent forms of this conjecture have been extensively studied
 by Kirchberg (\cite{Ki}), and subsequently in \cite{Pi}, \cite {Ra1}, 
\cite{Ra2}, \cite{Br}, \cite{Oz}.

 In this paper we prove that the Connes embedding conjecture is equivalent to
 a statement on the structure of the trace inequalities on matrices. To prove this, we deduce an analogue of the Hilbert $17$-th problem
in the context of type II$_1$ von Neumann algebras.

 More precisely, let $Y_1, \ldots,Y_n$ be $n$ indeterminates. 
 Let $\cI_n$ be the index set of all monomials in the variables
$Y_1, Y_2, \ldots, Y_n$
 \begin{equation*} \cI _n=\{(i_1, \ldots, i_p)\mid p \in \N, \, i_1, \ldots, i_p\in \{1, 2, \ldots, n\}\}\end{equation*}
 (we assume that $\emptyset \in \cI_n$ and that $\emptyset$ corresponds to the monomial which is identically one). For $I=(i_1, \ldots i_p)$ let $|\cI|=p$ and let $Y_I$ denote the monomial $Y_{i_1} Y_{i_2} \cdots Y_{i_p}$. For such an $I$ define $I^{\rm op} = (i_p, \ldots, i_1)$ and define an adjoint operation on $\C[Y_1, Y_2, \ldots, Y_n]$ by putting $Y^*_I = Y_{I^{\rm op}}$.

 We let $\can [Y_1, \ldots, Y_n]$ be the ring of all series
 \begin{equation*}
 V= \Big\{\sum_{I\in I_n} a_I Y_I\mid a_I\in \C, \, \Big \| \sum a_I Y_I \Big\|_R = \sum |a_I| R^{|I|}<\infty, \, \forall \, R>0\Big\}.\end{equation*} 
 It turns out (Section 2) that $V$ is a Fr\' echet space, and hence that $V$
 has a natural weak topology $\sigma (V, V^*)$. We will say that an element
 $q$ in $V$ is a sum of squares if $q$ is in the weak closure of the cone of
 sum of squares
 \begin{equation*}
 \sum_s p_s^* p_s, \quad p_s \in \can [Y_1, Y_2,\ldots, Y_n]. \end{equation*} 
 By $\can^{\rm sym} [Y_1, Y_2, \ldots, Y_n]$ we denote the real subspace of all analytic series that are auto-adjoint. 

 We say that two series $p, q$ in $\can^{\rm sym}[Y_1, Y_2, \ldots, Y_n]$ are {\it cyclic equivalent} if $p-q$ is a weak limit of sums of scalar multiples of monomials of the form $Y_I- Y_{\tilde I}$, where $I\in \cI_n$, and $\widetilde I$ is a cyclic permutation of $I$. 

 Our analogue of the Hilbert's $17$-th problem is the following 

 \begin{thm} Let $p \in \can ^{\rm sym}[Y_1, \ldots, Y_n]$ such that, whenever
   $M$ is a separable type {\rm II}$_1$ von Neumann algebra with 
faithful trace $\tau$ and $X_1, \ldots, X_n$ are selfadjoint elements in $M$, 
then by substituting $X_1, \ldots, X_n$ for $Y_1, \ldots, Y_n$, we obtain
$\tau (p(X_1, \ldots, X_n))\geq 0$. Then $p$ is cyclic equivalent to a weak limit of a sum of squares in $\can [Y_1, \ldots, Y_n]$.
 \end{thm}

 As a corollary, we obtain the following statement which describes the Connes
 embedding conjecture strictly in terms of (finite) matrix algebras. (A similar statement has been noted in \cite{DH}).

 \begin{cor} 
The Connes embedding conjecture holds if and only if whenever
$p\in \can^{\rm sym}[Y_1, \ldots, Y_n]$ is such that, whenever we substitute
selfadjoint matrices $X_1, X_2, \ldots, X_n$ in $M_N (\C)$ endowed with the
canonical  trace $\tr$, we have $\tr(p(X_1, X_2, \ldots, X_n))\geq 0$. Then $p$ should be equivalent to a weak limit of sums of squares.
 \end{cor} 

 To prove the equivalence of the two statements we use the fact that the Connes
 conjecture is equivalent to show that the set of non-commutative moments of
 $n$ elements in a type II$_1$ factor can be approximated, in a suitable way, by
 moments of $n$ elements in a finite matrix algebra.

 \section{Properties of the vector space $\can[Y_1, Y_2,  \ldots, Y_n]$}

 We identify in this section the vector space $\can [Y_1, \ldots, Y_n]$ with the vector space $V= \Big \{ (a_I)_{I\in \cI_n}\mid \sum\limits_{I\in \cI_n} |a_I| R^I< \infty , \, R>0\Big\}$. Denote by $\|(a_I)_{I\in \cI_n} \|_R = \sum \limits_{I \in \cI_n} |a_i| R^I$ for  $R>0$. Clearly, $\|\cdot\|_R$ is a norm on $V$. 
 In the next proposition (which is probably known to specialists, but we could not provide a reference) we prove that $V$ is a Fr\' echet space. For two elements $J, K \in \cI_n$, $J=(j_1, \ldots, j_n)$, $K= (k_1, \ldots, k_s)$, we denote by  $K\diez J= (j_1, \ldots, j_n, k_1, \ldots, k_s)$ the concatenation of $J$ and $K$.  
 \begin{prop}
 With the norms $\| \cdot \|_R$, $V$ becomes a Fr\'echet space. Moreover, the operation $*$ defined for $a= (a_I)$, $b= (b_I)$ by 
 \begin{equation*}
 (a* b)_I = \sum _{{\scriptstyle J,K \in \cI_n} \atop {\scriptstyle J\diez K=I}} a_Ja_K\end{equation*}
 (which in terms of $\can[Y_1, \ldots, Y_n]$ corresponds to the product of series) is continuous with respect to the norms $\|\cdot \|_R$, that is
 \begin{equation*}
 \|a * b\|_R\leq \|a\|_R \|b\|_R, \quad \mbox{ for all } R>0. \end{equation*}
 \end{prop}
 \begin{proof} To prove that $V$ is a Fr\'echet vector space, consider  a
   Cauchy sequence $b^s =(b_I^s)_{I\in \cI_n}$. Thus for all $R,\varepsilon>0$ there exist $N_{R, \varepsilon}$ in $\N$ such that for all $s, t> N_{R, \varepsilon}$ we have $\|b^s- b^t\|_R<\varepsilon$. Clearly, this implies that $(b^s_I)_{I\in \cI_n}$ converges pointwise to a sequence $(b_I)_{I\in \cI_n}$ and it remains to prove that $b$ belongs to $V$ and $b^s$ converges to $b$.

 To do this we may assume (by passing to a subsequence, and using a typical Cantor diagonalization process) 
 that 
 \begin{equation*} \|b^s- b^{s+1}\|_s \leq \frac 1 {2^s} \quad \mbox{for all $s$ in $\N$}.\end{equation*}

 But then for every $R>0$, we have that 
 \begin{align*}
 | b_I^s-b_I| R^{|I|} &\leq \sum_{{\scriptstyle s= [R] +1}\atop  {\scriptstyle I\in \cI_n}} |b^s_I-b^{s+1}_I| R^{|I|}\\&\leq \sum_{{\scriptstyle s=[R]+1}\atop {\scriptstyle I\in \cI_n}} | b_I^s - b_I^{s+1}| s^{|I|} \leq \sum _{s= [R]+1} \frac 1 {2^s}.\end{align*}
  This proves that $b\in V$ and that $b^s$ converges to $b$. 

  For the second part of the statement we have to estimate, for all $a, b\in V$ and $R>0$, the quantity
  \begin{align*}
  \sum_{I\in \cI_n} \sum_{{\scriptstyle J, K\in \cI_n}\atop {\scriptstyle J\diez K=I}} |a_J| \, |b_K| R^{|I|}&= \sum_{I\in \cI_n} \sum_{{\scriptstyle J,K \in \cI_n}\atop{\scriptstyle J\diez K=I}} |a_J|\, |b_K|R^{|J|+|K|}\\&=
  \bigg ( \sum_{J\in \cI_n} |a_J| R^{|J|}\bigg)\bigg( \sum _{K\in \cI_n}|b_K|R^{|K|}\bigg)= \| a\|_R \|b \|_R.\end{align*}
  This proves that $a*b\in V$ and that $\|a*b\|_R \leq \| a\|_R \|b\|_R$. \end{proof}

  In what follows we describe the dual $V^*$ of the Fr\'echet space $V$.
  \begin{lem} $V^*$ is identified with the space of all sequences of complex numbers \begin{equation*}\{(t_I)_{I\in cI_n} \mid  \exists\, R>0, \, \sup\limits _{I\in \cI_n}|t_I|R^{-|I|} < \infty\}.\end{equation*} 
  The duality between $V$ and $V^*$ is realized via the pairing 
  \begin{equation*}   \langle a, t\rangle = \sum_{I\in \cI_n} a_I t_I\quad
    \mbox{for } a\in V,\, t\in V^*,\end{equation*}
  which is convergent if $\sup\limits _{I\in \cI_n} |t_I|R^{-|I|}<\infty$ for some $R$.
  \end{lem}
  \begin{proof} It is obvious that each sequence $(t_I)_{I\in \cI_n}$ such that $\sup\limits_{I\in \cI_n} |t_I|R^{-|I|}$ defines an element in $V^*$.

  Conversely, if $\varphi$ is a continuous linear functional on $V$, then there exists a semi-norm $\|\cdot\|_R$ and a positive constant $C>0$, such that
  \begin{equation*} |\varphi ((a_I)_{I\in \cI_n})| \leq C \sum_{I \in \cI_n} |a_I|R^{|I|}.\end{equation*}
  But then by the usual duality between $\ell^1$ and $\ell^\infty$ and 
since $\{ (a_I)R^{|I|})_{I\in \cI_n}\mid (a_I)_{I \in \cI_n} \in V\}$ 
is dense in $\ell^1(\cI_n)$, it follows that there exists 
$(t_I)_{I\in \cI}$ such that $\sup \limits_{I\in \cI_n} (t_I) R^{-|I|}<\infty$
and such that $\varphi ((a_I)_{I\in \cI_n}) = \sum\limits_{I\in \cI_n}
a_It_I$.
\end{proof}

  By $V_{\rm sym}$ we consider the real subspace of $V$ consisting of all
  $\{(a_I)_{I\in \cI_n}\in V\mid a_{I^{\rm op}}= \overline {a_I} , 
\, \forall I\in \cI_n\}$, and by  $V^*_{\rm sym}$ we consider the space of all real functionals on $V_{\rm sym}$. 

  If $\sigma$ is a permutation of $\{1,2, \ldots, p\}$ (which we denote by $\sigma \in S_p$) and $I\in \cI_n$, $I=(i_1, i_2, \ldots, i_p)$ then by $\sigma(I)$ we denote $(i_{\sigma(1)}, i_{\sigma(2)}, \ldots, i_{\sigma(p)})$. By $S_{p,{\rm cyc}}$ 
  we denote the cyclic permutations as $S_p$. We omit the index $p$ when it is obvious.

Let $W = \{(t_I)_{I\in \cI_n} \mid t_I = t_{\sigma(I)},\,
\forall\, \sigma\in S_{\rm cyc} and \,\exists \, R>0$ such that
$\sup |t_I|R^{-|I|} < \infty, \, t_I = \overline{t_{I^{\rm op}}},
\, \forall\, I \in \cI_n\}$.

 \begin{lem} 
$W$ is a closed subspace of $V^*_{\rm sym}$. Moreover, if $\varphi _s=(t_I^s)_{I\in \cI_n} \in W$ converges to $\varphi \in W$ in the $\sigma (V^*, V)$ topology, then there exists $R>0$ such that 
 \begin{equation*} \sup_{s, I\in \cI_n}|t_I^s| R^{-|I|} <\infty.\end{equation*}
 \end{lem}
 \begin{proof} $W$ is a closed subspace of $V^*$ follows immediately from the fact that $W$ is the fixed point subspace for the actions of the finite groups $S_{p, {\rm cyc}}$, $p\geq 1$, on the components of the indices. The same is true for the third condition in the definition of $W$.

 The second statement is a standard consequence of the Banach-Steinhaus theorem. Indeed, if $\varphi_s \to \varphi$ in the $\sigma(V^*, V)$ topology, then $(\varphi _s)_{s\in \N}$ forms an equicontinuous family and hence there exist a semi-norm $\|\cdot\|_R$ on $V$ and a constant $C>0$ such that
 \begin{equation*} |\varphi_s((a_I)_{I\in \cI_n})| \leq C\| (a_I)_{I\in \cI_n}\|_R.\end{equation*}

   By the preceding section this means exactly the condition in the statement.
   \end{proof}

 We note the following straightforward consequence of the Bipolar theorem (\cite{Ru}, \cite{Simai-Reed}):

 \begin{lem} Let $K_m\subseteq K_2\subseteq W$ be closed convex cones, (which will correspond to the sets 
$K_{\rm mat}$ and $K_{{\rm II}_1}$ in the next section).  Let $L_2=K^0_2$, $L_m=K^0_m$ be
the corresponding polar sets (the polars are with respect $V_{\rm sym}$) in $V_{\rm sym}$. (In particular $W^\perp \subseteq L_2\subseteq
L_m$.)

 By the annihilator $W^\prep$ of the space $W$ we mean the (relative) annihilator with respect  to$V_{\rm sym}$. Likewise we do for the polar sets.

 Let $L_p$ be a $\sigma (V_{\rm sym}, V^*_{\rm sym})$ closed conex subset of $L_2$. Then to prove that $L_p+W^\perp=L_2$ is equivalent to prove that $L^0_p\cap W=K_2$ and hence it is sufficient to verify that $L^0_p \cap W\subseteq K_2$.
 \end{lem} 

 \begin{proof}
 Indeed, in this situation $L_p+W^\perp \subseteq L_2$ and since $L^0\cap W=
 (L_p+W^\perp)^0$ it follows immediately from the Bipolar theorem that
 $K_2=L_2^0\subseteq (L_p + W^\perp)^0 =L_p\cap W$. Moreover, if the last two
 closed convex sets are equal, i.e.\ if $K_2=L_p^0 \cap W$, then by the Bipolar theorem we get $L_2=K_2^0= L_p ^0 \cap W= (L_p+W^\prep )^{00}= L_p+W^\perp.$\end{proof}


\section{The set of non-commutative moments for variables \\ in a type II$_1$
von Neumann algebra}

We consider the following subsets of the (real) vector space $W$ introduced
in the preceding section:

Let 
\begin{align*}
K_{{\rm II}_1}  = \sigma(V^*,V)\mbox{-{\rm closure}}& \{(\tau(x_I))_{I\in \cI_n}
\mid M  \mbox{ type ${\rm II}_1$ separable von Neumann algebra, }\\
& \mbox{ with faithful trace } \tau, \; x_i = x_i^* \in M, \, i = 1,2,\ldots,n\}.
\end{align*}
Here by $x_\phi$ we mean the unit of $M$.

We also consider the following set. Let $\tr$ be the normalized trace
$\frac{1}{N}\trr$ on $M_N(\C)$.
Then 
$$
K_{\rm mat} = \sigma(V^*,V)\mbox{-{\rm closure}}\{(\tr(x_I))_{I\in \cI_n} \mid
N \in \N, \, x_i = x_i^* \in M_N(\C)\}.
$$
We clearly have that $K_{\rm mat} \subseteq K_{{\rm II}_1} \subseteq W$ and
as proven in \cite{Ra1}, $K_{\rm mat}$, $K_{{\rm II}_1}$ are convex 
$(\sigma(V^*,V)\mbox{-{\rm closed}}$ or $\sigma(V_{\rm sym}^*,V_{\rm sym})$ 
closed sets. Let $\tilde K_{{\rm II}_1}$, $\tilde K_{\rm mat}$ be the (convex) cones generated
by the convex sets $K_{{\rm II}_1}$ and $ K_{\rm mat}$

We consider also the following subset of $V_{\rm sym} \subseteq V$
$$
L_p = \overline{\rm co}^{\sigma(V,V^*)} \{a^* * a \mid a \in V\},
$$
where $*$ is the operation introduced in Proposition 2.1.

\begin{obs}
In the identification of $V$ with $\C_{\rm an}[Y_1,Y_2,\ldots,Y_n]$,
$L_p$ corresponds to the series that are limits of sums of squares. Moreover,
$W^\perp \cap V_{\rm sym}$ corresponds to the $\sigma(V,V^*)\mbox{-{\rm closure}}$ of the 
span of the selfadjoint part of 
monomials
of the form $Y_I - Y_{\sigma(I)}$, $\sigma \in S_p$. 
Hence $L_p + (W^\perp \cap V_{\rm sym}$ corresponds to the series in $\C_{\rm an}^{\rm sym}
[Y_1,\ldots,Y_n]$ that are equivalent to a weak $(\sigma(V,V^*))$
limit of sum of squares: $$\sum\limits_{s\in S} b_{s^*}^*b_s,\quad 
b_s \in \C_{\rm an}[Y_1,Y_2,\ldots,Y_n].$$

Moreover, $L_p \subseteq K_{{\rm II}_1}^0 \subseteq K_{\rm mat}^0$.
\end{obs}

\begin{proof} Indeed the series corresponding to $a^* * a$ is 
\begin{align*}
\sum_{I\in \cI_n} (a^* * a)_I Y_I & = \sum_I 
\sum_{J\sharp K} a_J^* a_K Y_{J\sharp K}  =
\sum_{I\in \cI_n} \sum_{J\sharp K=I} \overline{a_{J^{\rm op}}} a_K Y_J Y_K 
\\ & =
\sum_{J,K\in \cI_n} \overline{a}_{J} a_K Y_{J^{\rm op}} Y_K  =
 \bigg[ \sum_{J\in \cI_n} (\overline{a}_J Y_J^*)\bigg]^*
\bigg( \sum_{K\in \cI_n} a_K Y_K\bigg).
\end{align*}

Moreover, when pairing a series $\sum a_I Y_I \in 
\C_{\rm an}[Y_1,\ldots, Y_n]$
with $t_I = \tau(x_I)$,
$I \in \cI_n$, $x_1,\ldots,x_n \in M$, $x_i = x_i^*$, the result is
$\sum a_I \tau(x_I) = \tau\Big(\sum a_I x_I\Big)$. Hence
an element in $L$ paired with an element in $K_{{\rm II}_1}$ 
is the trace of a sum squares
and hence is positive. Thus, $L_p \subseteq K_{{\rm II}_1}^0$ and hence 
$L_p + (W^\perp \cap V_{\rm sym})
\subseteq K_{{\rm II}_1}^0$ so $L_p^0 \cap W \supseteq K_{{\rm II}_1}$.
\end{proof}

With these observations we can prove Theorem 1.1.

\begin{proof} By the above observations
and by Lemma 2.4 this amounts to prove that
$$
L_p + (W^\perp \cap V_{\rm sym}) = K_{{\rm II}_1}^0
$$
and by the previous remarks we only have to prove that
$$
L_p^0 \cap W \subseteq \tilde K_{{\rm II}_1}.
$$
Thus we want to show that if an element $\theta=(\theta_I)_{I \in \cI}$ with $\theta_{\emptyset}=1$ belongs
to $L_p^0$ then there exists a type II$_1$ von Neumann algebra 
$M$ with
faithful trace $\tau$, and selfadjoint elements $x_1,\ldots,x_n$ in $M$
such that
$$
\theta_I = \tau(x_I) \quad \mbox{ for } I \in \cI_n.
$$

To construct the Hilbert space on which $M$ acts consider
the vector space 
$$
\cH_0 = \C_{\rm an}[Z_1,\ldots,Z_n]
$$
where $Z_1,Z_2,\ldots,Z_n$ are indetermined variables. We consider the following
scalar product on $\cH_0$:
$$
\langle Z_I,Z_J \rangle = \theta_{J^{\rm op}\sharp I} 
$$
or more general for $(a_I)_{I \in \cI_n}$, $(b_I)_{I \in \cI_n}\in V$
$$
\bigg\langle \sum_{I \in \cI_n} a_I Z_I, 
\sum_{J \in \cI_n} b_J Z_J \bigg \rangle = \langle b^* * a, \theta \rangle.
$$
(Recall that $\langle \;,\;\rangle$ is the pairing between 
$V$ and $V^*$ introduced in Section 2).
Note that the scalar product is positive since
$$
\bigg\langle \sum_{I \in \cI_n} a_I Z_I, 
\sum_{I \in \cI_n} a_I Z_I \bigg \rangle = \langle a^* * a, \theta \rangle,
$$
which is positive since $t\in L_p^0 \cap W\subseteq L_p^0$.

We let $\cH$ to be the Hilbert space completion of $\cH_0$ after we mod
out by the elements $\xi$ with $\langle \xi,\xi\rangle =0$.
Note that in this Hilbert space completion,  we have obviously that
$\sum\limits_{I \in \cI_n} a_I Z_I$ is the Hilbert space limit after
$p \to \infty$ of $\sum\limits_{I \in \cI_n\atop |I|\leq p} a_I Z_I$.

We consider the following unitary operators $U_t^i$, $i=1,2,\ldots,n$, $t\in \R$
acting on $\cH_0$ isometrically and hence on $\cH$.

For $i \in \{1,2,\ldots,n\}$ and $q \in \N$ we will denote 
by $i^q$ the element $(i,i,\ldots,i)$ of length $q$ and belonging to
$\cI_n$.

The following formula defines $U_t^i$
$$
U_t^i Z_I = \sum_{s=0}^\infty \frac{(it)^s}{s!} Z_{i^s \sharp I}
\quad \mbox{ for } I \in \cI_n.
$$
To check that $U_t^i$ is corectly defined, and extendable 
to an unitary in $\cH$ it is sufficient to check that
$$
\langle U_t^i Z_I, U_t^i Z_J \rangle = \langle Z_I, Z_J \rangle
\eqno(1)
$$
for all $I,J \in \cI_n$. But we have
\begin{align*}
\langle U_t^i Z_I, U_t^i Z_J \rangle 
& = \sum_{p,q}\frac{(it)^p \overline{(it)^q}}{p!q!} 
\langle Z_{i^p\sharp I}, Z_{i^q \sharp J} \rangle 
= \sum_{p,q}\frac{(it)^p \overline{(it)^q}}{p!q!} 
\theta_{J^{\rm op} \sharp i^{p+q}\sharp I} \\
& = \sum_{p,q}\frac{i^{p-q}t^{p+q}}{p!q!} 
\theta_{J^{\rm op} \sharp i^{p+q}\sharp I}.
\end{align*}
Note that all the changes in summation are allowed 
and that all the series are absolutely
convergent since $(t_I)_{I\in \cI} \in W$ and hence there exists
$R>0$ such that $|t_I| \leq R^{|I|}$ 
(we also use the analiticy of the series
$\sum\limits_{p\geq 0} \frac{(it)^p}{p!}$.

For a fixed $k$, the coefficient of $t^k$ in the above sum is
$$
t^k(\theta_{J^{\rm op} \sharp i^k \sharp I}) 
\sum_{p+q=k}\frac{i^{p-q}}{p!q!} $$
and this vanishes unless $k=0$ because of the corresponding property of the
exponential function: $e^{it} e^{-it}=1$.

This proves (1) and hence that $U_t^i$ can be extended to an isometry on $\cH_0$
and hence on~$\cH$.

Next we check that for all $j=1,2,\ldots,n$,
$$
U_t^j U_s^j = U_{t+s}^j \quad \mbox{ for all } t,s \geq 0
$$
for $I \in \cI_n$. But
$$
U_t^j U_s^j Z_I = \sum_{p,q=0}^\infty\frac{(it)^p(is)^q}{p!q!} Z_{j^p\sharp j^q \sharp I}
= \sum_{k=0}^\infty \sum_{p+q=k} 
\frac{(it)^p(is)^q}{p!q!} Z_{j^k\sharp I}
$$
and this is then equal to $U_{t+s}^j Z_I$ because of the corresponding property
for the exponential function: 
$e^{is} e^{it}=e^{i(s+t)}$.
From this it then follows that $U_t^j$ is a one parameter group of unitaries for all
$j=1,2,\ldots,n$.

We let $M$ be the von Neumann algebra of $B(\cH)$ generated by all
$U_t^j$, $j=1,2,\ldots,n$, $t\in \R$. We will verify that $Z_\emptyset$ 
is a cyclic
trace vector for $M$ (and hence it will be also separating, \cite{Sz}).

To verify that $Z_\emptyset$ is a trace vector it sufficient to check that for all
$p,q$ and all $i_1,i_2,\ldots,$ $i_p, j_1,j_2,\ldots,j_q  \in \{1,2,\ldots,n\}$
we have that for all $t_1,\ldots,t_p, s_1,\ldots,s_q  \in \R$
$$
\langle U_{t_1}^{i_1} \cdots U_{t_p}^{i_p} \cdots 
U_{s_1}^{j_1} \cdots U_{s_q}^{j_q} Z_\Phi, Z_\Phi \rangle =
\langle U_{s_1}^{j_1} \cdots U_{s_q}^{j_q} U_{t_1}^{i_1} \cdots U_{t_p}^{i_p} 
Z_\Phi, Z_\Phi \rangle.\eqno(2)
$$
It is immediate that the right side of (2) is equal to 
$$
\sum_{\alpha_1,\ldots,\alpha_p =1 \atop \beta_1,\ldots,\beta_q =1}^\infty
\frac{(it_1)^{\alpha_1}\cdots (it_p)^{\alpha_p} (is_1)^{\beta_1}\cdots (is_q)^{\beta_q}}
{\alpha_1!\cdots\alpha_p! \beta_1! \cdots \beta_q!}
\theta_{i_1^{\alpha_1}\cdots i_p^{\alpha_p} j_1^{\beta_1}\cdots j_q^{\beta_q}}.
$$
Because of the cyclic symmetry of $(\theta_I)_{I\in \cI_n}$ is then equal to the lefthand
side of (2).

The sums involved are absolutely convergent since $|t_I| \leq R^{|I|}$, 
$I\in \cI_n$ for some\break $R>0$ and because of the fact that the scalar
function involved are entire analytic functions.

We now prove that $Z_\Phi$ is a cyclic vector. Denote by $iy_j$ the selfadjoint
generator for the group $U_t^j$, $j = 1,2,\ldots,n$.

We claim that for all $I\in \cI_n$, the vector $Z_I$ belongs to the domain
of $iy_j$ and that $y_j Z_I = Z_{j \sharp I}$, $j = 1,2,\ldots,n$,
$I\in \cI_n$. Indeed this follows by evaluating the limit 
$$
\lim_{t \to 0} \frac{1}{i} \frac{1}{t}(U_t^j Z_I - Z_I) - Z_{j \sharp I}.
$$
But this is equal to
$$
t\bigg( \sum_{q=2}^\infty \frac{(it)^q}{q!} Z_{j^n \sharp I}\bigg)
$$
and this converges to zero when $t\to 0$ because of the summability condition
for $(t_I)_{I\in \cI_n}$.

Indeed the square of norm of this term is 
$$
t^2 \sum_{p,q=2}^\infty \frac{(it)^p \overline{(it)^q}}{p!q!} 
\theta_{I^* \circ j^{n+m} \circ I}.
$$
From this we deduce that $Z_\Phi$ belongs to domain of $y_{j_1}\cdots y_{j_q}$
and $y_{j_1}\cdots y_{j_q}Z_\Phi = Z_{(j_1\cdots j_q)}$.

This implies in particular that $Z_\Phi$ is cyclic for $M$ and also that
$$
\langle y_I Z_\Phi, Z_\Phi \rangle = \theta_I, \quad I\in \cI_n.
$$
By \cite{Sz} it follows that $Z_\Phi$ is also separating for $M$ and hence
that $M$ is a type II$_1$ von Neumann algebra and that the vector form
$\langle \,\cdot\, Z_\Phi,Z_\Phi \rangle$ is a trace on $M$.

By definition the elements $y_j$, $j = 1,2,\ldots,n$ are selfadjoint and
affiliated with $M$.

Moreover since $\tau(m)=\langle m Z_\Phi,Z_\Phi \rangle$ is
a faithfull finite trace on $M$ and since there exist $R>0$ such that
$$
\tau(y_j^{2k}) = \theta_{j^{2k}} \leq R^{2k} \quad \mbox { for all } k,
$$
it follows that $y_j$ are bounded elements in $M$ such that
$\tau(y_I) = \theta_I$ for all $I\in \cI_n$.
This proves that $(\theta_I)_{I\in \cI_n}$ belongs to $K_{{\rm II}_1}$ and hence
that $L_p + (W^\perp \cap V_{\rm sym}) 
\subseteq \tilde K_{{\rm II}_1}$. This completes the proof of the
theorem.
\end{proof}

We are now ready to prove Corollary 1.2.

\begin{proof} With our previous notations this amounts to show that
the statement $L_p + (W^\perp \cap V_{\rm sym}) =  K_{\rm mat}^0$ is equivalent to the Connes
embedding conjecture.

Since we already know that $L_p + (W^\perp \cap V_{\rm sym})
\subseteq \tilde K_{{\rm II}_1}$, this amounts
to prove that the statement that $K_{\rm mat} = K_{{\rm II}_1}$ is equivalent
to Connes conjecture.

The only non-trivial part is that the statement $K_{\rm mat} = K_{{\rm II}_1}$ 
implies the Connes conjecture. So assume that for a given type II$_1$ von 
Neumann algebra $M$ and for given $x_1,\ldots,x_n$ in $M$ there exists
$(X_i^k)_{i=1}^n$ in $M_{N_k}(\C)$ such that 
$(\tr(X_I^k))_{I\in \cI_n}$ converges  in the $\sigma(W,V)$ topology to
$(\tau(x_I))_{I\in \cI_n}$.

By Lemma 2.3, we know that there  exists $R>0$ such that for all $k \in \N$ 
$$
|\tr(X_I^k)|\leq R^{|I|}
$$
and hence that
$$
\tr((X_j^k)^{2p}) \leq R^{2p} \quad \mbox{ for all } j = 1,2,\ldots,n,\,
p,k\in \N.
$$
But this implies that $\|X_j^k\|\leq R$ for all $k\in \N$, $j = 1,2,\ldots,n$.

But then the map $x_j  \to (X_j^k)_{k\in \N}$
extends to an algebra morphism from $M$ into $R^\omega$ which preserves the trace.
\end{proof}

\section*{Acknowledgements}
We are indebted to Professors R. Curto, F.~Goodman, C.~Frohman, and
C. Procesi
for several disscutions during the elaboration of this paper.
The author is grateful to the EPSRC and to the Deparment
of Mathematics of the Cardiff University, where part of this
work was done.


\begin{thebibliography}{aa}

\bibitem{Br}
{\sc N.P. Brown}, Connes' embedding problem and Lance's WEP, 
{\it Int. Math. Res. Not.} 2004, no. 10, 501--510.

\bibitem{Co}
{\sc A. Connes}, Classification of injective factors. Cases $II\sb{1},$
$II\sb{\infty },$ $III\sb{\lambda },$ $\lambda \not=1$, 
{\it Ann. of Math. (2)}
{\bf 104}(1976), no.~1, 73--115.

\bibitem{DH}
{\sc D. Hadwin},
A noncommutative moment problem,  {\it Proc. Amer. Math. Soc.} {\bf 129} (2001), no. ~6, 1785--1791


\bibitem{Ki}
{\sc E. Kirchberg}, 
Discrete groups with Kazhdan's property ${\rm T}$ and
factorization property are residually finite, {\it Math. Ann.} 
{\bf 299}(1994), no.~3, 551--563.

\bibitem{McD}
{\sc ?. McDuff}, Dusa Central sequences and the
hyperfinite factor, 
{\it Proc. London Math. Soc. (3)} {\bf 21}(1970), 443--461.

\bibitem{Oz}
{\sc N. Ozawa}, About the QWEP conjecture,
preprint, math.~OA/0306067.

\bibitem{Pi}
{\sc G. Pisier}, 
A simple proof of a theorem of Kirchberg and related results on 
$C\sp *$-norms,
{\it J. Operator Theory} {\bf 35}(1996), no.~2, 317--335.

\bibitem{Procesi}
{\sc C. Procesi}, The invariant theory of $n\times n$ matrices, 
{\it Advances in Math.} {\bf 19}(1976), no.~3, 306--381

\bibitem{Ra1}
{\sc F. R\u adulescu},
Convex sets associated with von Neumann algebras and Connes' approximate
embedding problem,
{\it Math. Res. Lett.} {\bf 6}(1999), no. 2, 229--236.

\bibitem{Ra2}
{\sc F. R\u adulescu},
The von Neumann algebra of the non-residually finite Baumslag group
$\langle a,b \mid a b^3 a^-1 = b^2 \rangle$ embeds into $R^\omega$
preprint, math.~OA/0004172.

\bibitem{Ru}
{\sc W. Rudin}, {\it Functional Analysis}, Second edition, 
International Series
in Pure and Applied Mathematics, McGraw-Hill, Inc., New York, 1991.

\bibitem{Simai-Reed} 
{\sc B. Simon, M. Reed},
Methods of modern mathematical physics. II. Fourier analysis, self-adjointness, Academic Press,
New York-London, 1975. 


\bibitem{Sz}
{\sc S. Str\u atil\u a, L. Zsid\'o},
Lectures on von Neumann algebras, 
Revision of the 1975 original, translated from the Romanian by
Silviu Teleman, Editura Academiei, Bucharest; Abacus Press, Tunbridge
Wells, 1979.

\end{thebibliography}
\end{document}